\documentclass{article}
\usepackage{epigraph}
\usepackage{listings}
\usepackage{amssymb}
\begin{document}


\title{Automatic classification of automorphisms of lower-dimensional  Lie algebras}
\author{C. Wafo Soh$^{1,2}$\\
$^1$Department of Mathematics  and Statistical Sciences,\\ Jackson State University\\
 JSU Box 17610, 1400 JR Lynch Street,\\ Jackson, MS 39217, USA\\ 
$^{2}$DST-NRF Centre of Excellence in \\Mathematical and Statistical Sciences,\\ School of Computer Science and Applied Mathematics,\\ University of the Witwatersrand, \\Johannesburg, Wits 2050, South Africa}
\maketitle
\begin{abstract}
We implement two algorithms  in MATHEMATICA for classifying automorphisms of lower-dimensional non-commutative Lie algebras. The first algorithm is a brute-force approach whereas the second is an evolutionary strategy. These algorithms are delivered as the MATHEMATICA package {\tt cwsAutoClass}. In order to facilitate the application of this package to symmetry Lie algebras of differential equations, we also provide a package, {\tt cwsLieSymTools}, for manipulating finite-dimensional Lie algebras of vector fields. In particular, this package allows the computations of Lie brackets, structure constants, and the visualization of commutator tables. Several examples are provided to illustrate the pertinence of our approach.   	
\end{abstract}


\section{Motivation}	
A symmetry of an object is a transformation that does not affect it. In several applications, it is desirable to compute all the symmetries of the object of interest. However, this is a difficult and sometimes 
an intractable problem. In the case of differential equations, Sophus Lie \cite{lie1891} discovered that the knowledge of its continuous symmetries unravels an algorithmic path for its integration by quadratures. Lie derived an algorithm for the computation of these symmetries, which in most cases, involved solving an overdetermined (i.e. there are more equations than unknowns) system of linear partial differential equations (PDEs).Until recently, there was not an algorithm for computing all the discrete symmetries of differential equations except  for the direct method which leads to equations as difficult to solve as the initial  equation itself. That is, the linearization afforded by continuous symmetries is lost in the search for discrete symmetries. Beside, the direct method for the search of discrete symmetries does not guarantee that one has found all the discrete symmetries unless one solves the nonlinear determining equations for them. Hydon \cite{hyd2000} discovered that, when the continuous symmetry Lie algebra of a differential equation is non-Abelian, in the symmetry group, discrete and continuous symmetries are entangled. However,  if one wants to compute all its discrete symmetries, one must first disentangle its  two types of symmetries  through  the classification of automorphisms of its  continuous symmetry Lie algebra. This disentanglement produces conditional determining equations for discrete symmetries that are first-order  quasilinear PDEs. This conditional determining  system is overterdermined and in principle can be solved using Lagrange-Charpit method of characteristics.       
Once, its solutions are obtained, the discrete symmetries are found by constraining these solutions to leave invariant the underlying differential equation. This last stage is a mere application of the chain rule for differentiation. However, there are three main roadblocks in the implementation of Hydon's algorithm for the computation of discrete symmetries of a differential equation admitting a non-Abelian symmetry Lie algebras. They are: (1) the solution of a large overtedermined  system of quadratic equations, (2) the classification of the automorphisms of the symmetry Lie algebra, and  (3) the solution of the conditional determining system. 
In this work, we focus primarily on the classification of automorphisms of lower-dimensional non-Abelian Lie algebras. Specifically, We implement in MATHEMATICA a brute-force approach together  with an evolutionary strategy for classifying the automorphisms of a non-commutative finite-dimensional Lie algebra.

We have structured this paper as follows. There are four sections including this introduction. In Section 2, we state the problem to be solved while fixing notations. Section 3 deals with the implementation of two automorphisms classification algorithms {\it viz.}  brute-force and genetic strategies. Section 4 is concerned with the test of our implementation. There, we consider several examples that illustrate diverse aspects of our code. The last section i.e. Section 5, summarizes our work and engages in some discussions about it.

\section{Statement of the problem}
Our goal in this section is to present the algorithm for classifying the automorphisms of a finite-dimensional non-Abelian  Lie algebra. We start with some preliminaries on Lie algebras \cite{graaf1997}.  

A Lie algebra $L$ over a field $\mathbb{F}\in\{\mathbb{R},\mathbb{C}\} $ is a vector space which is equipped with a bilinear operation $[.,.]$ (Lie bracket) which enjoys the following two properties: 
(1) for all $x,\,y \in L$, $[x,\,y]=-[y,\ x]$ (antisymmetry), and  
(2) for all $x,\, y$, and $z$, $[[x,\,y],\,z] + [[y,\,z],\,x] +[[z,x],y] =0$ (Jacobi's identity). The dimension of a Lie algebra is its dimension when it is treated as a vector space. In the sequel, we shall be dealing solely with finite-dimensional Lie algebras. So suppose that $L$ is  such Lie algebra with basis $\{e_1,e_2,\ldots,e_n  \}$. For all $i,j = 1:n$, $[e_i, e_j] \in L$. So, for all $i,j = 1:n$, $[e_i, e_j] = \sum_{k=1}^nc_{ij}^k e_k$ for some constants
$c_{ij}^k \in \mathbb{F},\, k= 1:n$. The $c_{ij}^k$'s are known as the {\em structure constants} of  $L$. It can be verified that, thanks to the antisymmetry of the Lie bracket and Jacobi's identity, the structure constants have the following properties: for $i,j, k, l = 1:n$,
\begin{equation}
c_{ij}^k = - c_{ji}^k \mbox{ and } 
\sum_{m=1}^n (c_{ij}^m c_{mk}^\ell +c_{jk}^m c_{mi}^\ell +c_{ki}^m c_{mj}^\ell)=0. \label{sp1}
\end{equation} 

A linear map, $f$,  from the Lie algebra  $L$ to itself is called and automorphism if  it is a bijection and  for all $x,y\in L$, $f([x,\,y]) = [f(x),\, f(y)]$. If we designate by $B=[b^j_i]$ the matrix of $f$, where $i$ and $j$ are respectively the row and column numbers, then the entries of $B$ are constrained by the two conditions:
\begin{eqnarray}
\det(B) &\ne & 0, \label{sp2}\\
\sum_{\ell=1}^{n} c_{ij}^{\ell}b^{k}_{\ell} & = &
\sum_{\ell,m =1}^{n} c^k_{\ell m} b^{\ell}_i b^m_j,\quad
i<j =1:n.\; k = 1:n, \label{sp3}
\end{eqnarray}
Conversely, it can be shown that if a matrix $B$ satisfies the constraints (\ref{sp2})-(\ref{sp3}), then it induces an automorphism of the Lie algebra $L$. In the system of equations (\ref{sp3}), there are $n^2(n-1)/2$ equations for $n^2$ unknown. As $n$ becomes large, it can be a daunting task to solve it. Beside, one may encounter an explosion of cases which also erode considerably our computational budget. Once the system (\ref{sp3}) is solved, one must filter the solution set according to the constraint (\ref{sp2}).

For all $j =1:n$, the mapping $x\mapsto [x, e_j]$ is a derivation  of $L$ whose matrix $C(j)$ is such that its entry at position $(i,k)$ is $(C(j))^{k}_i = c_{ij}^k$. Each such derivation with matrix $C(j)$, $j=1:n$, induces a  one-parameter family of automorphisms of $L$ with matrices 
$\{ \exp(\epsilon C(j)) : \epsilon \in \mathbb{F} \}$. Let us introduce the notation $A(j,\epsilon) =  \exp(\epsilon C(j))$ for latter convenience. Also, we shall denote the set of all the matrices $B$ satisfying Eqs. (\ref{sp2})-(\ref{sp3}) by $\mbox{Aut}(L)$. Thus, for all $j=1:n$ and $\epsilon \in \mathbb{F}$, $A(j,\epsilon)\in \mbox{Aut}(L)$. 

We define on $\mbox{Aut}(L)$ a relation $\sim_d$ by: $B_1 \sim_d B_2$ if and only if 
there are $(\epsilon_1, \epsilon_2, \ldots, \epsilon_n )\in \mathbb{F}^n$ and a permutation $\sigma$ of $\{1,2,\ldots, n\}$ such that 
\begin{equation}
B_1 = A(\sigma(1),\epsilon_1)\, A(\sigma(2), \epsilon_2)\,\ldots A(\sigma(n), \epsilon_n)\, B_2.  \label{sp4}
\end{equation} 
It can be shown that the relation $\sim_d$ is reflexive and symmetric i.e. it is a {\em dependency} relation. Our main goal is to determine the dependency classes modulo $\sim_d$. That is, we want to find all the {\em traces} of the dependency $\sim_d$. For a given trace, we wish to select its sparsest member as representative. Note that there are three types of $A(j, \epsilon)'s$: (1) those that are equal to the identity matrix $I_n$, (2) those that are diagonal matrices distinct from the identity, and (3) those that are non-diagonal. The $A(j,\epsilon)$'s that are identity matrices come from the $C(j)$'s that are zero matrices. They are generated by $e_j$'s that commute with all the elements of $L$ i.e. the $e_j's$ that belong to the center of $L$. They do not affect the definition of $\sim_d$. So they may be disregarded. Thus, we are left with the nontrivial diagonal and nondiagonal $A(j, \epsilon)$'s. Since a diagonal matrix commute with any matrix of the appropriate size, in the definition of  $\sim_d$, after discarding the identities, we may reposition all the diagonal matrices in front and their order does not matter. We may assume without lost of generality that the diagonal $A(j, \epsilon)$'s are numbered from $1$ to $s$, and the nondiagonal ones are numbered from  $s+1$ to $p$, with $p\le n$. Now, we may rephrase the definition of $\sim_d$ in the following way: $B_1 \sim_d B_2$ if and only if there are $(\epsilon_1,\epsilon_2,\ldots, \epsilon_p)\in \mathbb{F}^n$  and a permutation $\tau$ of $\{s+1, s+2,\ldots, p\}$ such that
\begin{equation}
B_1 = A(1,\epsilon_1)\,\ldots A(s, \epsilon_s)
A(\tau(s+1),\epsilon_{s+1})\ldots A(\tau(p),\epsilon_p) \, B_2.  \label{sp5}
\end{equation}
Given $B_2\in \mbox{Aut}(L)$, we are after its sparsest  dependents $B_1$ which contain less parameters (i.e unspecified entries). We accomplish this by selecting an appropriate permutation $\tau$ and by choosing  $\epsilon_{s+1}$ to $\epsilon_{p}$ such that as much as possible entries of  $B_3=A(\tau(s+1),\epsilon_{s+1})\ldots A(\tau(p),\epsilon_p) \, B_2$ are equal to zero. Then, we pick $\epsilon_1$ to $\epsilon_s$  such that some rows or columns of   $A(1,\epsilon_1)\,\ldots A(s, \epsilon_s)B_3 = B_1$ are appropriately scaled. 

For a fixed $B \in \mbox{Aut}(L)$ and a  permutation $\tau$ of $\{s+1, s+2, \ldots, p\}$, we may end up with several sparse dependents. We denote by $\tau(B)$ the set of such dependents. We define the {\em fitness} of $\tau(B)$ as the geometric mean of the sparsity of its elements. As $\tau$ ranges over permutations of $\{s+1, s+2, \ldots, p\}$, we select the permutation which provides the fittest set of dependents of $B$.

In the scheme for finding dependency classes modulo $\sim_d$ that  we have just described, there are  two main roadblocks: (1) the solution of the system (\ref{sp2})-(\ref{sp3}), and (2) the search of optimal representative of dependency classes. Indeed the complexity of both these problems is exponential in the dimension, $n$, of $L$. Thus as $n$ increases, it is imperative to adopt strategies for taming this issue. For the first problem, we may partially solve the system (\ref{sp2})-(\ref{sp3}) by relying upon the structure of its equation: We may start by solving linear  and quadratic equations with fewer terms. 
 As for the second problem, we adopt an evolutionary strategy at the expenses of having sub-optimal representatives of traces.

In the remainder of this paper, we work through our implementation of the computation of dependency classes modulo $\sim_d$, test our implementation by treating few examples, and  discuss the virtues and limitations of our implementation while suggesting possible future improvements.        
\section{Implementation}
This section reifies the algorithm described in the previous section using MATHEMATICA. We assume that the reader is familiar with MATHEMATICA's syntax. We adopt mostly a functional programming approach. In all the code snippets provided, build-in functions start with a capital letters whereas our variables start with lower cases. We will not provide helper functions since they can be consulted in the accompanying source code.   

\subsection{Solution of the system (\ref{sp2})-(\ref{sp3})}
The function that solves the system (\ref{sp2}) -(\ref{sp3}) is 
named {\tt simplifiedB}. It consumes five inputs in the following order: (1) the dimension of the Lie algebra $L$, {\tt dim}, (2) a function {\tt c} that allows the calculation of structure, (3)  a symbol {\tt b} that is used to name the entries of the matrix $B$, (4) the field, {\tt dom}, of the Lie algebra $L$, and (4) a list of exigences, {\tt constraints}, on the entries of the matrix $B$. It is coded using some helper functions that we do not provide.
\begin{lstlisting}[numbers = left, breaklines=true]  
simplifiedB[dim_,c_,b_, dom_:Complexes,constraints_:{}]:=simplifiedB[dim,c,b, dom,constraints]= Module[{ sol = solveSys[dim,c,constraints][b], 
bb = symbolicB[dim,b], bs}, bs=  Map[bb/.#&,sol];
If[dom === Reals, Select[bs, realMatrixQ], bs]]//Simplify;
\end{lstlisting}

\subsection{Determination of a single dependency class}
Given a solution {\tt b} of the system (\ref{sp2})-(\ref{sp3}) and a permutation {\tt sigma} of  the nondiagonal $A(j,\epsilon)$'s of the Lie algebra $L$, the function {\tt reducedB} finds all its sparsest dependents. The inputs {\tt c, dim} and {\tt dom } are as before. The main challenge in implementing the function {\tt reducedB} is to be able deal with all the cases. Indeed as we zero an entry of $B$ by multiplying it by an $A(j, \epsilon)$, we may end up with several possible $\epsilon$'s which must be treated separately. The same remark applies when we scale a row or column of $B$ using a non-trivial diagonal $A(j, \epsilon)$.

\begin{lstlisting}[numbers =left, breaklines=true]
reduceB[c_,dim_, perm_, dom_: Complexes,constraints_:{}][b_] :=
reduceB[c,dim, perm,dom,constraints][b]=
Module[{ f, g, auto = gatherByDiagonal[allStructMatrix[dim, c] ],k, sb={b}},
$Assumptions = Fold[And, Flatten[Table[Element[b[i,j], dom],{i,1,dim}, {j,1,dim}]]];
auto[[1]] = Select[auto[[1]],  # != zeroMatrix[dim]&]; 
auto[[2]] = Permute[auto[[2]],perm]; 
f = Join@@Map[findEpsAndReduceB[k,0,dom,constraints] , myOuter[Dot, #1, MatrixExp[k*#2] ]]&;
g = Join@@Map[findEpsAndReduceB[k,1,dom ,constraints] , myOuter[Dot, #1, MatrixExp[k*#2]]]&;  
{Fold[ g, Fold[f, sb, auto[[2]]],auto[[1]]], perm }]//Simplify       
\end{lstlisting}
\subsection{Determining all the dependency classes}
Now that we know how to compute  the dependents of a solution $B$ of the system (\ref{sp2})-(\ref{sp3}) for a given permutation of the non-diagonal $A(j, \epsilon)$'s, we can select the fittest when this permutation ranges over the set of all such permutations. As we mentioned earlier, the fitness of a trace is taken as the geometric mean of the sparsities of its members.
\subsubsection{The brute-force approach}
\epigraph{When in doubt, use brute force.}{Ken Thompson}
 The function that compute the optimal set of dependents for all the solutions of (\ref{sp2})-(\ref{sp3}) is called {\tt autoClassificationBruteForce}. Its inputs are similar to those of previous functions except for the new input {\tt name}  which is a string that will be used as radical for renaming  some variable entries of dependents. The outputs is a pair  whose first entry is the set of representatives of dependency classes whereas the second entry comprises restrictions encountered during computations. We employed the latter output mainly for backtracking calculations.

\begin{lstlisting}[numbers = left, breaklines = true]
autoClassificationBruteForce[dim_ ,c_, dom_: Complexes, constraints_:{}][b_, name_:\[Alpha]]:= Module[{s, cond, f},
f=  Not[Normal[Det[#]]===0]&;
s = Map[autoClassificationBruteForceOneB[dim,c,dom,constraints], simplifiedB[dim,c,b,dom,constraints]];
cond = Fold[Or, DeleteDuplicates@Map[#[[2]]&,s]];
s = DeleteDuplicates[Join@@Map[renameVariable[name, Flatten[symbolicB[dim,b]],#,dom]&, Join@@Map[First, s]]];
{Select[s, f] ,  cond}];  
\end{lstlisting} 

\subsubsection {An evolutionary approach}
\epigraph{
I have called this principle, by which each slight variation, if useful, is preserved, by the term of Natural Selection.
}{Charles Darwin}

In the function {\tt autoClassificationBruteForce}, as the number of nondiagonal inner automorphisms increases, our computational budget is quickly depleted. In oder to palliate this situation, we adopt an evolutionary strategy through the function {\tt autoClassificationGen}. Besides the previous type of inputs, it consumes the following ones: {\tt sigma, p, popSize, numGen} and {\tt cond}. The argument {\tt sigma} is the percentage of the optimal fitness we want to achieve, and {\tt p} is the mutation probability. The initial population size is {\tt popSize} whereas {\tt numGen} is the largest  number of generations one is willing to go through before stopping if the  desired fitness is not realized. In the code below, the function 
{\tt aiReduceB} does the heavy lifting of our genetic strategy which  encompasses our mating, mutation and selection schemes. 
For a quick intuition into our evolutionary scheme, consider two parents $par_1 =(dep_1, \mu_1)$ and $par_2 = (dep_2, \mu_2)$ selected at random  with  probabilities proportional to their respective fitnesses, where the $dep_i$'s are lists of dependents and the $\mu_i$'s are the corresponding permutations of non-diagonal inner automorphisms (the so-called $A(j,\epsilon)$'s). In each individual's definition pair, we shall refer to the first entry as its phenotype and the  second entry as its DNA. Theses biological analogies are self-explanatory.  Note that, given  an individual begotten from 
a solution of (\ref{sp2})-(\ref{sp3}) and its DNA, we can always recover its phenotype through computations.  Now, the parents 
$par_1$ and $par_2$ produce six  possible types of offspring with DNAs $\mu_1^{-1},\, \mu_2^{-1}, \mu_1 \circ \mu_2, \,
\mu_2\circ \mu_1,\, \mu_1 \circ \mu_2 \circ \mu_1^{-1}$, and
$ \mu_2 \circ \mu_1 \circ \mu_2^{-1}$. After possible mutation of children, the fittest amongst parents and offspring earns the right to belong to the next generation. We accomplish mutation by simply swapping two randomly chosen entries of the underlying DNA. In MATHEMATICA syntax, a permutation is represented by a list. For instance a permutation of ${1,2,3}$ may be represented as $\mu =\{3, 1, 2 \}$. It means that $\mu(1) =3$, $\mu(2) =1$, and $\mu(2) =2$. Thus, $\mu' = \{2,1,3\}$ is possible mutation of $\mu$ which is produce by swapping the nucleotides  $3$ and $2$.  

The output of the function {\tt autoClassificationGen} is formatted according to the boolean input {\tt cond}. It is a list of three or two elements according to whether {\tt cond} is true or false. When {\tt cond} is true, the output comprises the list of dependence classes, the constraints encounter during computations as well as the DNAs of dependency classes. In the even {\tt cond} is false, the second entry in the output list is omitted.
 
\begin{lstlisting}[numbers = left, breaklines = true]
autoClassificationGen[sigma_, p_, c_, dim_, popSize_:5, numGen_:1, dom_:Complexes, cond_: False,constraints_:{}][b_,name_:\[Alpha]]:= Module[{s,perm,cd,f},
f = Not[Normal[Det[#]]=== 0]&;
s= Map[aiReduceOneB[sigma,p,c,dim,popSize, numGen, dom,constraints],  simplifiedB[dim,c,b,dom,constraints]];  
perm = DeleteDuplicates[Map[#[[2]]&,s]]; 
s  = Join@@Map[  Join[#[[1]]]&,s];
s = Map[collectConditions, s];              
{s,cd} = unZip[s]; 
s = DeleteDuplicates[Join@@Map[renameVariable[name, Flatten[symbolicB[dim,b]],#, dom]&,s]]; 
cd = Fold[Or,False, Map[Simplify, DeleteDuplicates[Flatten[ Map[Fold[And ,True, #]&, cd]]]]]; 
If[cond, {Select[s,f], cd, perm}, {Select[s,f], perm}]];  
\end{lstlisting} 
\section{Tests}
Here we test our implementation of the classification of automorphisms of non-commutative finite-dimensional Lie algebras. We have encapsulated our algorithms in a MATHEMATICA package called {\tt cwsAutoClass}. Additionally, we provide a package {\tt cwsLieSymTools} which facilitates some computations pertaining to finite-dimensional Lie algebra represented in terms of vector fields. In particular it allows the calculations of the structure constants and the generation of the commutator table given a basis of the  underlying Lie algebra.
The examples we shall treat are done in MATHEMATICA  version 11 run on a DELL INSPIRON laptop with WINDOWS 10 operating system and the following additional specifications: INTEL CORE i3-3227U @ 1.90 GHz processor, 3.96 GB of usable RAM  and a 64-bit operating system. 

\subsection{Dependency classes of 3D non-Abelian Lie Algebras}
There are ten non-commutative  three-dimensional Lie algebras \cite{pawi}. We provide in this section the commands for computing their dependency classes using the package {\tt cwsAutoClass}.   

In first line of the code provided below, replace the comment with the appropriate information before compilation. The instruction of that line loads the package {\tt cwsAutoClass}.

\begin{lstlisting}[numbers = left, breaklines = true]
Get[(* Put the location of the file cwsAutoClass.m here e.g. "C:\\Users\\Celestin\\Desktop\\Trip_to_SA\\cwsAutoClass.m" *)]; 

c1[1, 2, 2] = 1; c1[2, 1, 2] = -1; c1[i_, j_, k_] := 0;
Map[MatrixForm, 
autoClassificationBruteForce[3, c1, Reals][b, \[Theta]][[1]]] // 
Quiet // Timing

c2[2, 3, 1] = 1; c2[3, 2, 1] = -1; c2[i_, j_, k_] := 0;
Map[MatrixForm, 
autoClassificationBruteForce[3, c2, Reals][b, \[Theta]][[1]]] // 
Quiet // Timing

c3[1, 3, 1] = 1; c3[3, 1, 1] = -1; c3[2, 3, 1] = 1; c3[3, 2, 1] = -1; 
c3[2, 3, 2] = 1; c3[3, 2, 1] = -1; c3[i_, j_, k_] := 0;
Map[MatrixForm, 
autoClassificationBruteForce[3, c3, Reals][b, \[Theta]][[1]]] // 
Quiet // Timing

c4[1, 3, 1] = 1; c4[3, 1, 1] = -1; c4[2, 3, 2] = 1; c4[3, 2, 2] = -1; 
c4[i_, j_, k_] := 0;
Map[MatrixForm, 
autoClassificationBruteForce[3, c4, Reals][b, \[Theta]][[1]]] // 
Quiet // Timing

c5[1, 3, 1] = 1; c5[3, 1, 1] = -1; c5[2, 3, 2] = -1; c5[3, 2, 2] = 1; 
c5[i_, j_, k_] := 0;
Map[MatrixForm, 
autoClassificationBruteForce[3, c5, Reals][b, \[Theta]][[1]]] // 
Quiet // Timing

Clear[a]; $Assumptions = (-1 < a  < 1) &&  a != 0;  c6[1, 3, 1] = 1; 
c6[3, 1, 1] = -1; c6[2, 3, 2] = a; c6[3, 2, 2] = -a; 
c6[i_, j_, k_] := 0;
Map[MatrixForm, 
autoClassificationBruteForce[3, c6, Reals][b, \[Theta]][[1]]] // 
Quiet // Timing

c7[1, 3, 2] = -1; c7[3, 1, 2] = 1; c7[2, 3, 1] = 1; c7[3, 2, 1] = -1; 
c7[i_, j_, k_] := 0;
Map[MatrixForm, 
autoClassificationBruteForce[3, c7, Reals][b, \[Theta]][[1]]] // 
Quiet // Timing

Clear[a];  $Assumptions = a > 0; c8[1, 3, 1] = a; c8[3, 1, 1] = -a; 
c8[1, 3, 2] = -1; c8[3, 1, 2] = 1; c8[2, 3, 1] = 1; c8[3, 2, 1] = -1; 
c8[2, 3, 2] = a; c8[3, 2, 2] = -a; c8[i_, j_, k_] := 0;
Map[MatrixForm, 
autoClassificationBruteForce[3, c8, Reals][b, \[Theta]][[1]]] // 
Quiet // Timing

c9[1, 2, 1] = 1; c9[2, 1, 1] = -1;  c9[2, 3, 3]  = 1; 
c9[3, 2, 3] = -1; c9[3, 1, 2] = 2;  c9[1, 3, 2] = -2; 
c9[i_, j_, k_] := 0;
Map[MatrixForm, 
autoClassificationBruteForce[3, c9, Reals][b, \[Theta]][[1]]] // 
Quiet // Timing

c10[1, 2, 3] = 1; c10[2, 1, 3] = -1;  c10[3, 1, 2]  = 1; 
c10[1, 3, 2] = -1; c10[2, 3, 1] = 1;  c10[3, 2, 1] = -1; 
c10[i_, j_, k_] := 0;
lst = Map[MatrixForm, 
autoClassificationBruteForce[3, c10, Reals][b, \[Theta]][[1]]] //
Quiet // Normal // Timing
\end{lstlisting}
\subsection{Dependency classes of some 4D non-Abelian Lie algebras}
The code provided below computes the dependency classes of the non-decomposable non-commutative Lie algebras $A_{4,1}$, $A_{4,4}$, $A_{4,9}^1$, and $A_{4,12}$ \cite{pawi}.
\begin{lstlisting}[numbers = left, breaklines = true]
Get[(* Put the location of the file cwsAutoClass.m here e.g. "C:\\Users\\Celestin\\Desktop\\Trip_to_SA\\cwsAutoClass.m" *)]; 

c1[2, 4, 1] = 1; c1[4, 2, 1] = -1; c1[3, 4, 2] = 1; 
c1[4, 3, 2] = -1;   c1[i_, j_, k_] := 0;
Map[MatrixForm, 
autoClassificationBruteForce[4, c1, Reals][b, \[Mu]][[1]]] // 
Quiet // Timing

c2[1, 4, 1] = 1; c2[4, 1, 1] = -1; c2[2, 4, 1] = 1; c2[4, 2, 1] = -1; 
c2[2, 4, 2] = 1; c2[4, 2, 2] = -1; c2[3, 4, 1] = 1; c2[4, 3, 1] = -1; 
c2[3, 4, 3] = 1; c2[4, 3, 3] = -1; c2[i_, j_, k_] := 0;
Map[MatrixForm, 
autoClassificationBruteForce[4, c2, Reals][b, \[Mu]][[1]]] // 
Quiet // Timing

c3[2, 3, 1] = 1; c3[3, 2, 1] = -1; c3[1, 4, 1] = 2; c3[4, 1, 1] = -2; 
c3[2, 4, 2] = 1; c3[4, 2, 2] = -1; c3[3, 4, 3] = 1; c3[4, 3, 3] = -1; 
c3[i_, j_, k_] := 0; 
Map[MatrixForm, 
autoClassificationBruteForce[4, c3, Reals][b, \[Mu]][[1]]] // 
Normal // Quiet // Timing

c4[1, 3, 1] = 1; c4[3, 1, 1] = -1; c4[2, 3, 2] = 1; c4[3, 2, 1] = -1; 
c4[1, 4, 2] = -1; c4[4, 1, 2] = 1; c4[2, 4, 1] = 1; c4[2, 4, 1] = -1; 
c4[i_, j_, k_] := 0; 
Map[MatrixForm, 
autoClassificationBruteForce[4, c4, Reals][b, \[Mu]][[1]]] // 
Normal // Quiet // Timing
\end{lstlisting}
\subsection{Dependency classes of symmetry Lie algebras}
Here, we treat the dependency classes of symmetry Lie algebras of various PDEs. We shall employ the package {\tt cwsLieSymTools} to facilitate the computation of structure constants and the visualization of commutator tables. We consider in turn the symmetry Lie algebras of the spherical Burgers \cite{doen1990} ,  Harry-Dym \cite{here1994},  and Black-Scholes \cite{ibra,silb2005}  equations. 
\begin{lstlisting}[numbers = left, breaklines = true]
Get[(* The location of cwsAutoClass goes here e.g. "C:\\Users\\Celestin\\Desktop\\Trip_to_SA\\cwsAutoClass.m"*)]; 
Get[(* The location of cwsLieSymTools.m goes here e.g. "C:\\Users\\Celestin\\Desktop\\Trip_to_SA\\cwsLieSymTools.m" *)];
(* Spherical Burgers' equation *)
vars = {t, x, u}; X1 = {-2*t, -x, u}; X2 = {0, Log[t], 1/t}; X3 = {0, 
1, 0}; listSym = {X1, X2, X3}; paramList = {}; dim = 3;
commutatorTable[vars, listSym, X, paramList, Reals] // Quiet
c1 = structureConstant[vars, listSym, paramList, Reals];
Map[MatrixForm, 
autoClassificationBruteForce[dim, c1, Reals][b, \[Alpha]][[1]]] // 
Quiet // Timing
Map[MatrixForm, 
autoClassificationGen[0.99, 0.0001, c1, dim, 150, 20, Reals, True][
b, \[Alpha]][[1]]] // Quiet // Timing

(* Harry-Dym equation *)
vars = {t, x, u}; H1 = {0, 1, 0}; H2 = {0, x, u}; H3 = {0, x^2, 
2*x*u}; H4 = {1, 0, 0}; H5 = {t, 0, -u/3}; listSym = {H1, H2, H3, 
H4, H5} ; paramList = {}; dim  = 5;
commutatorTable[vars, listSym, X, paramList, Reals] // Quiet
c2 = structureConstant[vars, listSym, paramList, Reals];
Map[MatrixForm, 
autoClassificationBruteForce[dim, c2, Reals][b, \[Alpha]][[1]]] // 
Quiet // Timing
Map[MatrixForm, 
autoClassificationGen[0.3, 0.0001, c2, dim, 150, 20, Reals, True][
b, \[Alpha]][[1]]] // Quiet // Timing

(* Black-Sholes equation *)
Clear @@ {A, d, c}; vars = {t, x, u}; Y1 = {1, 0, 0} ; Y2 = {0, x, 
0} ; Y3 = {2*t, (Log[x] + d*t)*x, 2*c*t*u} ; 
Y4 =   {0, A^2*t*x, (Log[x] - d*t)*u} ;   Y5 = {2*A^2*t^2, 
2*A^2*t*x*Log[x], ((Log[x] - d*t)^2 + 2*A^2*c*t^2 - A*t^2)*u};       
Y6 = {0, 0, 
u}; X1 = (1/A^2)*(Y1 + d*Y2 + 
c*Y6); X2 = Y2; X3 = Y3 ; X4 = Y4; X5 = (1/2) Y5; X6 = Y6; dim = 6;
listSym = {X1, X2, X3, X4, X5, X6}; paramList = {A, d, c};
commutatorTable[vars, listSym, X, paramList, Reals] // Quiet
c3 = structureConstant[vars, listSym, paramList, Reals];
Map[MatrixForm, 
autoClassificationBruteForce[dim, c3, Reals][b, \[Alpha]][[1]]] // 
Quiet // Timing
Map[MatrixForm, 
autoClassificationGen[0.99, 0.0001, c3, dim, 150, 20, Reals, 
True][b, \[Alpha]][[1]]] // Normal // Quiet // Timing
\end{lstlisting}

\section{Conclusion and discussions} 
We have implemented a MATHEMATICA package, {\tt cwsAutoClass}, that allows the automatic classification of automorphisms of non-Abelian lower-dimensional Lie algebras. It includes two functions {\tt autoClassification\-BruteForce} and {\tt autoClassificationGen} which respectively  implement brute force and evolutionary strategies. We demonstrated this package by considering several examples. Additionally, we have showcased a package, 
{\tt cwsLieSym\-Tools} which facilitates the calculation of structure constants and the visualization of commutator tables of finite-dimensional Lie algebras of vector fields. The latter package is particularly useful when discussing the automorphisms of symmetry Lie algebras.

The package {\tt cwsAutoClass} has some bottlenecks stemming from the solution of a large system of quadratic equations and  the exploration of a large search spaces to find fittest solutions. We address the second problem through a genetic strategy. However, the first issue remains. A possible way to tackle it consists in partially solving the system (\ref{sp2})-(\ref{sp3}) followed by our evolutionary algorithm. Then, solve the remaining equations of the system (\ref{sp2})-(\ref{sp3}) for each dependency class of the optimal solution. We shall in the future extend  our package {\tt cwsAtouClass}  along these lines.    
\section*{Acknowledgments}
I gratefully acknowledges partial  financial support form the National Research Foundation (NRF) of South Africa. I thank Prof. F M Mahomed for facilitating such funding. I shared portions of the code  of the package {\tt cwsLieSymTools } in 2016 with Ms. Nomsa Ledwaba who was then a Masters student at the University of Cape Town.

\end{document}